\documentclass[letterpaper, 10 pt, conference]{ieeeconf}
\IEEEoverridecommandlockouts
\usepackage{times}   
\usepackage{amsmath,amsfonts} 
\usepackage{amssymb} 
\usepackage{graphicx}
\usepackage{bm}
\usepackage{enumerate}
\usepackage{booktabs}
\usepackage[ruled,vlined,linesnumbered]{algorithm2e}

\newcounter{thmcount}

\newcounter{algcount}

\newtheorem{assump}[thmcount]{Assumption}

\newtheorem{alg}[algcount]{Algorithm}

\makeatletter
\newcommand{\removelatexerror}{\let\@latex@error\@gobble}
\makeatother

\allowdisplaybreaks[4]

\DeclareMathOperator*{\argmin}{argmin}

\def\N{\mathbb{N}}
\def\ones{{\bf 1}}
\def\ind{\mathcal{I}}
\def\lubx{[\underline{x}_k^{(i)},\,\bar{x}_k^{(i)}]}

\title{\LARGE \bf
Parallel ADMM for robust quadratic optimal resource allocation problems
}

\author{Zawar Qureshi$^{1}$ \and Sebastian East$^{1}$ \and Mark Cannon$^{1}$
\thanks{$^{1}$Department of Engineering Science, University of Oxford, Parks Road, Oxford OX1 3PJ, UK.}
\thanks{\texttt{qureshizawar@gmail.com}, 
\texttt{sebastian.east@eng.ox.ac.uk},
\texttt{mark.cannon@eng.ox.ac.uk}}%
}

\begin{document}
\maketitle

\pagestyle{empty}

\begin{abstract}
An alternating direction method of multipliers (ADMM) solver is described for optimal resource allocation problems with separable convex quadratic costs and constraints and linear coupling constraints. We describe a parallel implementation of the solver on a graphics processing unit (GPU) using a bespoke quartic function minimizer. An application to robust optimal energy management in hybrid electric vehicles is described, and the results of numerical simulations comparing the computation times of the parallel GPU implementation with those of an equivalent serial implementation are presented.   
\end{abstract}
\section{Introduction}

An important class of control problems can be formulated as optimal resource allocation problems with separable convex cost functions, separable convex constraints and linear coupling constraints. Examples include energy management in hybrid electric vehicles~\cite{koot05,buerger18b}, optimal power flow and economic dispatch problems~\cite{vandenbosch87}, 
portfolio selection~\cite{markowitz52}, 
and diverse problems in manufacturing, production planning and finance~\cite{patriksson08}.
%
The computation required to solve these problems grows with the length of the horizon over which performance is optimized, and the rate of growth becomes more pronounced for scenario programming formulations~\cite{campi08} that provide robustness to demand uncertainty. 
%
This limits the applicability of strategies such as model predictive control that require online solutions of resource allocation problems.

The alternating direction method of multipliers (ADMM) is a first order optimization method that can exploit separable cost and constraint functions, thus forming the basis of a solver suitable for parallel implementation~\cite{eckstein94,boyd11}.
This paper describes a graphics processing unit (GPU) parallel implementation of ADMM that avoids the growth in computation time with horizon that is observed in serial implementations. 
This is possible because GPUs have large numbers (typically thousands) of processing units and efficient memory architectures.
Parallel ADMM implemented on general purpose GPUs therefore provides a low-cost approach for time-critical optimal resource allocation problems.  

Under the assumption that the separable cost and constraint functions are defined by quadratic polynomials, we derive in this paper an ADMM iteration in which the most computationally demanding step is the minimization of a set of quartic polynomials. Iterative methods have previously been proposed for solving polynomial functions on the GPU~\cite{reis08,klopotek12}. However, to retain the convergence guarantees of ADMM, we devise an exact quartic minimization algorithm for GPUs. The solver was implemented on the GPU using CUDA~\cite{cuda_guide} and we compare its computation times with an equivalent serial implementation for a range of problem sizes. We also describe the application of ADMM to robust optimal energy management in a plugin hybrid electric vehicle (PHEV). The GPU implementation of this solver was tested in simulations against corresponding serial implementations. The results demonstrate a significant speed-up and indicate that the approach is
feasible for robust supervisory control of PHEVs with realistic horizon lengths.

\section{Problem Statement}
\label{sec:problem}

Consider a resource allocation problem in which a collection of resources, which could for example represent alternative power sources onboard a hybrid vehicle, is to be used to meet predicted demand
at each discrete time instant over a future horizon.
%
%
For a horizon of $n$ time steps, let $y_k$ denote the demand and $x^{(i)}_k$ denote the resource allocated from source $i$ at time step $k$, for $k \in\N_n:=\{1,\ldots,n\}$ and $i \in \N_m:=\{1,\ldots,m\}$. In order to meet the predicted demand, we require
\[
\sum_{i=1}^m x_k^{(i)} \geq y_k , \ \ k \in \N_n 
\]
where $x_k^{(i)}$ is subject to bounds, $x_k^{(i)}\in \lubx$, $k\in\N_n$ as well as constraints on total capacity:
\[
\sum_{k=1}^n g_k^{(i)}(x_k^{(i)}) \leq c^{(i)} , \ \ i \in \N_m .
\]
Capacity constraints represent, for example, the total energy, $\smash{c^{(i)}}$, available from the $i$th power source over the $n$-step horizon, with losses accounted for through the functions ${g_k^{(i)}}$. 

The cost associated with  $x_k^{(i)}$ is denoted $f_k^{(i)}(x_k^{(i)})$. This represents, for example, the operating cost or fuel cost (including losses) associated with the $i$th resource at time $k$. The objective is to minimize total cost subject to constraints:
\begin{equation}\label{eq:nominal_opt}
\begin{aligned}
\min_{x_k^{(i)}\in \lubx} & \ \
\sum_{i=1}^m \sum_{k=1}^n f^{(i)}_k(x_k^{(i)}) \\
\text{subject to}\quad \  & \ \
\begin{aligned}[t]
& \sum_{i= 1}^m x_k^{(i)}  \geq y_k   
&& k \in \N_n \\
& \sum_{k=1}^n g_k^{(i)}(x_k^{(i)}) \leq c^{(i)}
&& i \in \N_m . 
\end{aligned}
\end{aligned}
\end{equation}
This problem is convex if $f_k^{(i)}$ and $g_k^{(i)}$ are convex functions of their arguments. This is a common assumption in energy management problems with nonlinear losses~\cite{egardt14,hadj16}.

\begin{assump}\label{assump:convex}
The functions $f_k^{(i)}$ and $g_k^{(i)}$ are convex, for 
all $i \in\N_m$ and $k \in \N_n$.
\end{assump} 

If the demand sequence ${y = \{y_k\}_{k\in\N_n}}$ is stochastic, we assume that independent samples $y^{(j)}$, $j\in\N_q:=\{1,\ldots,q\}$, can be drawn from the probability distribution of $y$.

\begin{assump}\label{assump:sample}
The samples $y^{(j)} = \{y_k^{(j)}\}_{k\in\N_n,j\in\N_q}$ are such that $y^{(j)}$ and $y^{(l)}$ are independent for $l\neq j$.
\end{assump} 

Various robust formulations of the nominal allocation problem~(\ref{eq:nominal_opt}) are possible using a set of sampled sequences $y^{(j)}$, $j\in\N_q$. For example, an open loop formulation optimizes a single sequence $\{x_k^{(i)}\}_{k\in\N_n}$ for each $i\in\N_m$ with the predicted power demand sequence in~(\ref{eq:nominal_opt}) replaced by $\{\max_j y_k^{(j)}\}_{k\in\N_n}$. Alternatively, a receding (or shinking) horizon approach introduces feedback by solving the resource allocation problem at each sampling instant using current capacity constraints and demand estimates. In this case a sequence $\{x_k^{(i,j)}\}_{k\in\N_n}$ is computed for each $i\in\N_m$ and $j\in\N_q$ subject to the additional constraint that, for each $i$, the first element in the sequence should be equal for all $j$:
\begin{equation}\label{eq:robust_opt}%
\begin{aligned}
& \min_{x_1^{(i)}\in\mathbb{R}, \, x_k^{(i,j)}\in \lubx} \ \ 
\tfrac{1}{q}\sum_{j=1}^q\sum_{i=1}^m \sum_{k=1}^n f_k^{(i,j)} (x_k^{(i,j)}) \\
& \text{subject to} \ \ 
\begin{alignedat}[t]{2}
& \sum_{i= 1}^m x_k^{(i,j)}  \geq y_k^{(j)}  \ 
& k \in \N_n, \, j \in\N_q \\
& \sum_{k=1}^n g_k^{(i,j)}(x_k^{(i,j)}) \leq c^{(i)} \qquad
& i \in\N_m, \, j\in\N_q \\
& x_1^{(i,j)} = x_1^{(i)} \ 
& i\in\N_m, \, j\in\N_q
\end{alignedat}
\end{aligned}
\end{equation}
%
Although the objective in (\ref{eq:robust_opt}) is an empirical mean, the approach of this paper allows for alternative definitions such as the maximum  over $j\in\N_q$ of $\sum_{i,k} f_k^{(i,j)}(x_k^{(i,j)})$.
The probability that the solution of~(\ref{eq:robust_opt}) violates the constraints of~(\ref{eq:nominal_opt}) has a known probability distribution that depends in on the number, $q$ of samples and on the structural properties of (\ref{eq:robust_opt}) -- see e.g.~\cite{campi08,calafiore10} for details -- moreover
these results apply for any demand probability distribution. 

An ADMM solver for~(\ref{eq:robust_opt}) is summarised in Appendix~A. The approach exploits the separable nature of the cost and constraints of~(\ref{eq:robust_opt}), which are sums of functions that depend only on variables at individual time-steps, to derive an iteration in which the greatest contribution to the overall computational load is the minimization of scalar polynomials in (\ref{eq:admm_iter}a). Consequently the iteration does not require matrix operations (inversions or factorizations) and, as discussed below, it is suitable for implementation on extremely simple computational hardware.
For each $i\in\N_m$, the polynomial minimization in step~(\ref{eq:admm_iter}a) can be performed in parallel for ${k\in\N_n}$ and $j\in\N_q$, since the update for $\smash{x_k^{(i,j)}}$ does not depend on $\smash{x_{l}^{(i,r)}}$ for $(l,r)\neq (k,j)$. A parallel implementation therefore allows the ADMM iteration to be computed quickly even for long horizons and large numbers of scenarios. 

In the following sections, we describe a GPU-based parallel implementation of the ADMM iteration defined in (\ref{eq:admm_iter}a-i). We make the assumption that $\smash{f_k^{(i,j)}}$ and $\smash{g_k^{(i,j)}}$ are convex quadratic functions of their arguments. Under this assumption (\ref{eq:robust_opt}) becomes a quadratically constrained quadratic program, and the iteration (\ref{eq:admm_iter}a) then requires the minimization of a 4th order polynomial, for each $i,j$ and $k$.

\begin{assump}\label{assump:quadratic}
$f_k^{(i,j)}$ and $g_k^{(i,j)}$ are quadratic functions:
\begin{alignat*}{2}
f_k^{(i,j)}(x) &= \alpha_{k,2}^{(i,j)} x^2 + \alpha_{k,1}^{(i,j)} x + \alpha_{k,0}^{(i,j)} 
\\
g_k^{(i,j)}(x) &= \beta_{k,2}^{(i,j)} x^2 + \beta_{k,1}^{(i,j)} x + \beta_{k,0}^{(i,j)}
\end{alignat*}
for $i \in\N_m, \, j \in\N_q, \, k \in \N_n$.
\end{assump} 

\section{Parallelization}
\label{sec:parallelization}

This section discusses parallelization of the ADMM iteration in Appendix~A, focussing on the minimization step~(\ref{eq:admm_iter}a). Under Assumption~\ref{assump:quadratic}, this requires the minimization of 4th order polynomials, for which we describe a tailored algorithm in Section~\ref{sec:quartic}. 
We consider NVIDIA GPUs programmed in CUDA C/C++ \cite{cuda_guide} based on a single instruction multiple thread (SIMT) execution model. The CPU (host) and GPU (device) have separate memory, making it necessary to copy data between the host and device memory. Copying data back and forth adversely affects computation time, and the computational benefits of parallelization may therefore only become apparent for large problems. We therefore investigate in Section~\ref{sec:comp_quartic} the dependence of the computational speed-up that can be achieved via parallelization on the number of quartic polynomials to be minimized.

\subsection{CUDA programming model}

A CUDA program consists of code running sequentially on the host CPU and code that runs in parallel on the GPU.
A function (known as a  \textit{kernel}) involving $N$ parallel executions of an algorithm is implemented by executing $N$ separate CUDA threads. Each CUDA thread resides within a one-, two- or three-dimensional block of threads. 
%
A kernel can be executed by multiple equally-shaped thread blocks, which are themselves organized into a grid of thread blocks.
Each thread block is required to execute independently, concurrently or sequentially, and this creates a scalable programming model in which blocks can be scheduled on any available multiprocessor within a GPU. 
%
%
%
%
During execution, CUDA threads may access data from multiple memory spaces:
each thread has private local memory; each thread block has shared memory visible to all threads of the block and with the same lifetime as the block; and all threads have access to the same global device memory.
We use the runtime API (Application Programming Interface) to perform memory allocation and deallocation, and to transfer data between the host and device memory.
%
%

For the CUDA implementations discussed in Sections~\ref{sec:comp_quartic} and~\ref{sec:comp_admm}, individual kernels were devised for specific elementwise operations while cuBLAS library functions were used for vector-matrix multiplications. To compare  the performance of non-parallel implementations, serial algorithms were written in C++ with compiler optimizations (e.g.~vectorization where possible) set for speed (/Ox) unless stated otherwise. Intel Math Kernel Library BLAS routines were used for vector-matrix multiplications.
Both implementations were coded using Microsoft Visual Studio.

\subsection{Quartic minimization}
\label{sec:quartic}

Iterative solution methods, such as Newton's method~\cite{reis08} and recursive de Casteljau subdivision~\cite{klopotek12} have previously been proposed for solving systems of polynomials on a GPU. However, when used to minimize the quartic polynomial in~(\ref{eq:admm_iter}a), iterative methods can adversely affect ADMM convergence due to residual errors in their solutions. For parallel GPU implementations of the iteration~(\ref{eq:admm_iter}) therefore, a better approach is to compute the closed form solutions of the cubic polynomial obtained after differentiating the quartic in (\ref{eq:admm_iter}a) with respect to $\smash{x_k^{(i,j)}}$.
Algebraic solutions of cubic equations can be derived in various ways, including Cardano's method and Vieta's method~e.g.~\cite[Sec.~14.7]{dummit04}). We use a combination of these methods to first determine whether roots are real and/or repeated before they are calculated. 
For problems in which real solutions are expected, this allows significant computational savings and avoids the memory requirements of storing complex roots. 

Consider the quartic polynomial
\begin{equation}\label{eq:general_quartic_equation}
f(x) = Ax^4 + Bx^3 + Cx^2 + Dx +E 
\end{equation}
and let $b=3B/4A$, $c = C/2A$ and $d=D/4A$. The extremum points of $f$ are the solutions of the cubic equation 
\begin{equation}\label{eq:general_cubic_equation}
x^3 + bx^2 + cx +d = 0 .
\end{equation}
Let $Q$, $R$ and $\Delta$ be defined by
\begin{gather*}
Q = c/3 - b^2/9 , \qquad
R = bc/6 - b^3/27 - d/2  \\
\Delta = Q^3 + R^2 .
\end{gather*}
If $\Delta \leq 0$,  all three roots of~(\ref{eq:general_cubic_equation}) are real. In this case let
\[
\theta = \cos^{-1}\Bigl(R/\sqrt{-Q^3}\Bigr) ,
\]
then the roots are 
\begin{align*}
x_a &= 2\sqrt{-Q}\cos{\Bigl(\theta/3\Bigr)} - b/3\\
x_b &= 2\sqrt{-Q}\cos{\Bigl(\theta/3 +2\pi/3\Bigr)} - b / 3\\
x_c &= 2\sqrt{-Q}\cos{\Bigl(\theta/3 +4\pi/3 \Bigr)} - b / 3 .
\end{align*}
If $\Delta > 0$, there is only one real root. In this case, defining
\[
S = \sqrt[3]{R + \sqrt{\Delta}}, \qquad
T = \sqrt[3]{R - \sqrt{\Delta}} \, ,
\]
the real root is given by
\[
x_a = S + T - b / 3 .
\]
For the special case of $Q=R=0$, the three roots are real and repeated, and Vieta's formula gives the solution
\[
x_a = x_b = x_c = - b/3 .
\]
We additionally require the solver to determine which root of~(\ref{eq:general_cubic_equation}) is the minimizing argument of $f$. This is simplified by sorting the roots and discarding the middle root, which is necessarily a maximizer. The minimizer can then be determined by comparing the value of the objective function at the two remaining roots.
The quartic minimization algorithm 
is summarised in Algorithm~\ref{alg:cubic}.

\begin{figure}[htb]
\removelatexerror
\begin{algorithm}[H]
\SetKwInOut{Input}{Input}
\SetKwInOut{Output}{Output}
\SetKwRepeat{Do}{Do}{while}

\Input{Coefficients $A$, $B$, $C$, $D$}
\Output{Minimizer $x^\ast$}
Evaluate $b$, $c$, $d$, $Q$, $R$, $\Delta$\;
\uIf{$\Delta > 0$}{
Evaluate $S$, $T$\; 
$x^\ast \gets S+T - b/3$\;
}
\uElseIf{$Q=R=0$}{
$x^\ast \gets -b/3$\;
}
\Else{
Evaluate $x_a$, $x_b$, $x_c$\;
$(x_1,x_2,x_3) \gets\text{sort}(x_a,x_b,x_c)$\;
\mbox{$\delta\hspace{-0.1em} f \!\gets\! \frac{1}{4}(x_1^4 \!-\! x_3^4) \!+\! \frac{b}{3}(x_1^3\!-\!x_3^3) \!+\! \frac{c}{2}(x_1^2\! -\! x_3^2) \!+\! d(x_1\!-\!x_3)$}\;\!\!\!\!\!
\uIf{$\delta\hspace{-0.1em} f > 0$}{
$x^\ast \gets x_3$\;}
\uElse{
$x^\ast \gets x_1$\;}
}
\caption{Quartic minimizer}\label{alg:cubic}
\end{algorithm}
\end{figure}


\subsection{Comparison of CPU/GPU performance}
\label{sec:comp_quartic}

The results presented in this section and Section~\ref{sec:comp_admm} were obtained using an Intel i5-7300HQ CPU with base clock speed 2.5 GHZ, 12GB of RAM, and an Nvidia GTX 1060 3GB GPU. All operations were performed with single-precision floating point numbers. 

To test performance we generate $N$ quartic equations with randomly chosen coefficients. For the CPU implementation we solve the $N$ equations sequentially using Algorithm~\ref{alg:cubic}.
For the GPU implementation, all available threads (up to a limit of 18432 threads for the Nvidia GTX 1060 3GB) are scheduled to perform Algorithm~\ref{alg:cubic} in parallel.
The test procedure for the CPU-based algorithm was as follows.
\begin{enumerate}
  \item Generate $N$ random sets of quartic coefficients.
  \item Record start time $T_1$.
  \item Minimize the $N$ quartics using Algorithm~\ref{alg:cubic}. 
  \item Record end time $T_2$ and elapsed time $T$ = $T_2- T_1$.
\end{enumerate}

The test procedure for the GPU parallel cubic solver is similar to the CPU procedure with step 3 changed to:
\begin{enumerate}
  \setcounter{enumi}{2}
  \item
  \begin{enumerate}
  	\item Input data to the GPU.
  	\item Minimize the $N$ quartics using Algorithm~\ref{alg:cubic}
implemented on the GPU.
  	\item Copy results back to the CPU.
  \end{enumerate}
\end{enumerate}
The test procedure was repeated at least 10 times and the average value taken for the run-time. Tests were carried out with and without CPU compiler speed optimizations. 

\begin{figure}[t]
\centerline{\includegraphics[scale = 0.5, trim={8mm 0 17mm 8mm}, clip]{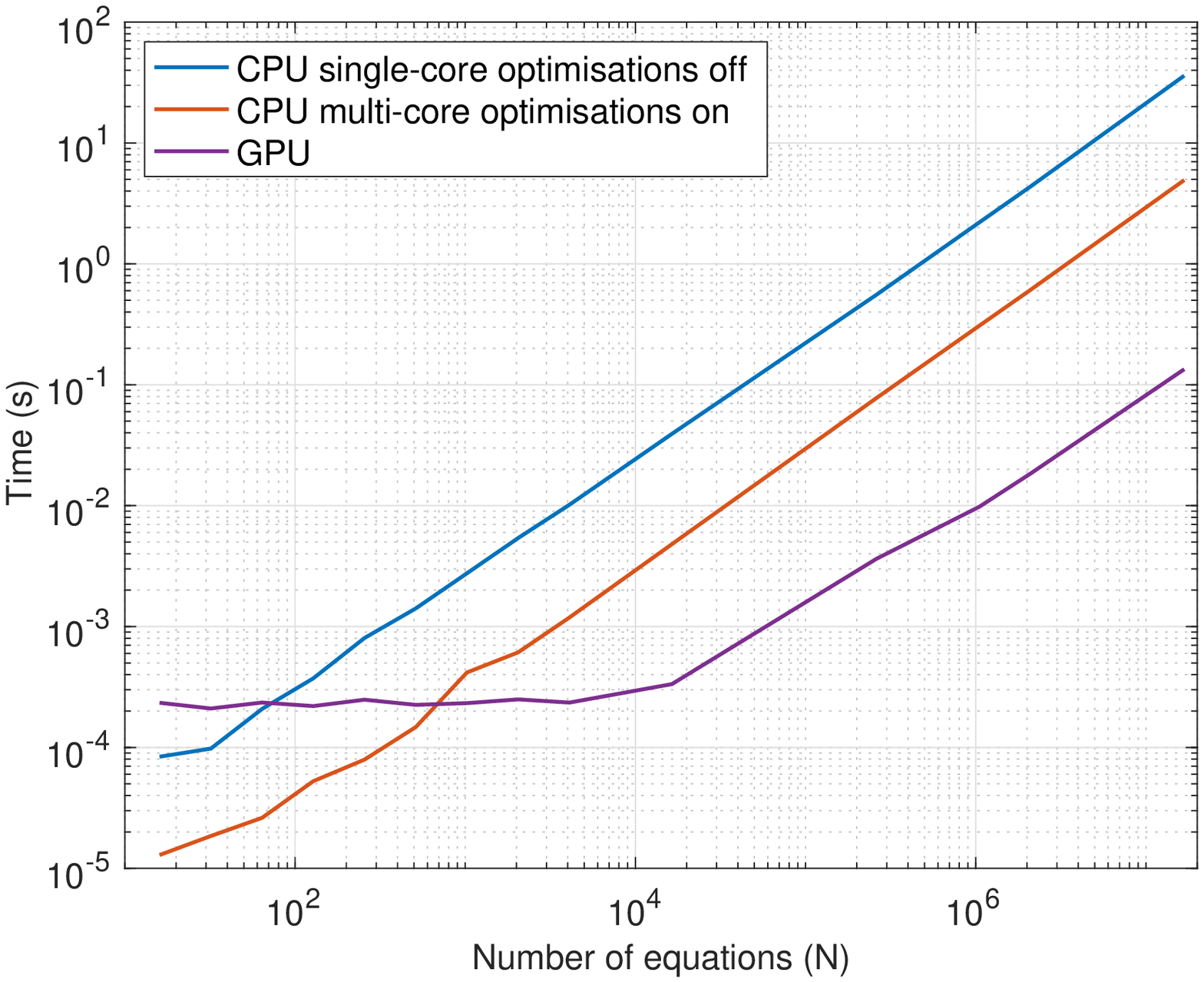}}
\vspace{-2mm}
\caption{Speed comparison for minimizing $N$ quartic equations. For CPU implementations, speed optimizations were either enabled or disabled using compiler flags (/Ox and /Od respectively).}
\label{fig:cubicBench}
\vspace{-2mm}
\end{figure}

Figure \ref{fig:cubicBench} compares the GPU and CPU execution speeds as $N$ varies. For $N > 10^6$, a $30$-fold increase in speed is obtained using the GPU in comparison with the CPU implementation with 4 cores. The speed-up factor rises to $210$ compared with a single-core CPU implementation with speed optimizations off. Although the CPU implementation is faster for $N < 10^3$, the GPU becomes significantly faster for large $N$. This is expected as the number of GPU threads must be sufficiently large to counter the effects of overheads such as memory transfer times. 


\section{PHEV Supervisory Control}
\label{sec:PHEVexample}

This section describes an application of the robust resource allocation problem~(\ref{eq:robust_opt}) and the parallelization approach of Section~\ref{sec:parallelization} to optimal energy management for a Plugin Hybrid Electric Vehicle (PHEV).
PHEVs employ hybrid powertrains that allow the use of more than one energy source. Here we consider a PHEV with two power sources: an electric motor/generator with battery storage, and an internal combustion engine.
The supervisory control problem involves determining the optimal power split between the internal combustion engine and the electric motor in order to  minimize a performance cost, defined here as the fuel consumption over a prediction horizon, while ensuring that the current and predicted future power demand is met and that battery stored energy remains within allowable limits.

We consider a parallel hybrid architecture~\cite{koot05,buerger18b} in which the propulsive power ($y$) demanded from the vehicle's powertrain is provided by the sum of power from the combustion engine ($x^{(1)}$) and power from the electric motor $(x^{(2)})$. Independent samples, $\smash{\{y_k^{(j)}\}_{j\in\N_q,k\in\N_n}}$, of the future demand sequence are assumed to be available, and the power balance
\[
y_k^{(j)} \leq x_k^{(1,j)} + x_k^{(2,j)}, \ \ j\in\N_q , \ \ k\in \N_n
\]
is imposed over an $n$-step future horizon.
%
For each demand scenario $j$ and time step $k$ the engine speed $\smash{\omega_k^{(j)}}$ is assumed to be known and equal to the electric motor speed whenever the clutch is engaged. 

The fuel energy consumed during the $k$th time interval is obtained from a quasi-static map, $\smash{f_k^{(1,j)}(x_k^{(1,j)})}$, relating engine power output to fuel consumption for the given engine speed $\smash{\omega_k^{(j)}}$.
Likewise, the battery energy output during the $k$th period, $\smash{g^{(2,j)}_k(x_k^{(2,j)})}$, is 
determined for a given $\omega_k^{(j)}$ by a quasi-static model of battery and electric motor losses.
The energy stored in the battery at the end of the horizon must be no less than a given terminal quantity, $E_n$, and hence the fall in battery stored energy over the prediction horizon must not exceed $\Delta E = E_0 - E_n$, where $E_0$ is the initial battery state.
%
%
No constraint is applied to total fuel usage and no cost associated with electrical energy usage
(i.e.~$g_k^{(1,j)} = 0$ and $\smash{f_k^{(2,j)}}=0$),
so the optimization problem~(\ref{eq:robust_opt}) becomes
\begin{equation}\label{eq:PHEV_opt}
\begin{aligned}
& \min_{x_1^{(i)}\in\mathbb{R},\  x_k^{(i,j)} \in\lubx}  \ \ 
 \tfrac{1}{q} \sum_{j=1}^{q}\sum_{k=1}^{n} f_k^{(1,j)}(x_k^{(1,j)}) 
\\
& \begin{alignedat}{3}
\text{subject to} \ \ 
& x_k^{(1,j)} + x_k^{(2,j)} \geq y_k^{(j)},  &\quad & j \in\N_q , \ k\in\N_n \\
& \sum_{k=1}^{n} g_k^{(2,j)}(x_k^{(2,j)}) \leq \Delta E, &\quad & j \in\N_q \\
& x_1^{(i,j)} = x_1^{(i)}, &\quad & i =1,2, \ j\in\N_q
\end{alignedat}
\end{aligned}
\end{equation}
The updates of $x_k^{(i,j)}$ in (\ref{eq:admm_iter}a) therefore require the minimization of a quadratic for $i=1$ and a quartic polynomial for $i=2$.
The bounds $\underline{x}_k^{(i)}$ are determined by constraints on power and torque output of powertrain components (see e.g.~\cite{buerger18b}); for example engine braking ($x^{(1,j)}_k < 0$) is ignored, so ${\underline{x}_k^{(1)} = 0}$. 
Assumption~\ref{assump:convex} (convexity) is satisfied if $f_k^{(1,j)}$ and $g_k^{(2,j)}$ are convex functions, which is usually the case for real powertrain loss maps (\cite{buerger18b}). 
%
In addition, cost and constraint functions are approximated by fitting quadratic functions to $\smash{f^{(1,j)}_k}$  and $\smash{g^{(1,j)}_k}$ at each time-step (see~\cite{buerger18b} for details), and Assumption~\ref{assump:quadratic} is therefore satisfied.

At each sampling instant a supervisory control algorithm based on the online solution of~(\ref{eq:PHEV_opt}) with a shrinking horizon can be summarised as follows.

\addtocounter{algcount}{1}
\begin{alg}[Supervisory PHEV control]\label{alg:PHEV_control}
\mbox{}\vspace{-0.5mm}
\begin{enumerate}
\item Draw samples, $\smash{\{y_k^{(j)}\}_{j\in\N_q,\, k\in\N_n}}$, from the distribution of the predicted power demand sequence.
\item For the current battery energy level, $E_0$, set $\Delta E = E_0 - E_n$, and compute for each $j\in\N_q$ the engine speed sequence $\{\omega_k^{(j)}\}_{k\in\N_n}$ corresponding to $\{y^{(j)}_k\}_{k\in\N_n}$, and hence determine $f_k^{(1,j)}$ and $g_k^{(2,j)}$ for $k \in\N_n$.
\item Solve~(\ref{eq:PHEV_opt}) by applying the ADMM iteration (\ref{eq:admm_iter}a-i) until the primal and dual residuals satisfy $r < \bar{r}$ and $\sigma < \bar{\sigma}$.
\item Apply the engine and electric motor power $x_1^{(1)}$, $x_1^{(2)}$.
\end{enumerate}
\end{alg}

\subsection{Modelling assumptions and optimization parameters}

We consider a $1900$\,kg vehicle with a hybrid powertrain consisting of a $100$\,kW gasoline internal combustion engine, $50$\,kW electric motor, and a $21.5$\,Ah lithium-ion battery.
Power demand sequences were generated by adding random sequences (Gaussian white noise of standard deviation $250$\,W and $50$\,rpm, low-pass filtered with cutoff frequency $0.02$\,Hz) to the power demand and engine speed computed for the FTP-75 drive cycle.
The sampling interval was fixed at~$1$\,s. 
The times at which the vehicle stopped were also modified to simulate variable traffic conditions.
This method of generating power demand is used simply to illustrate the computational aspects of the approach; more accurate methods of generating scenarios in this context (e.g.~\cite{jo17,ripaccioli10}) have been proposed in the literature on stochastic model predictive control.


The initial battery energy level $E_0$ was set at $60\%$ of battery capacity, $E_{\max}$, and the terminal level, $E_n$, at $50\%$ of capacity. Whenever the power demand $y^{(j)}_k$ is negative, $40\%$ is assumed to be recovered through regenerative braking. 
%

\begin{figure}[t]
\vspace{-2mm}
\centerline{\includegraphics[scale = 0.48, trim={5mm 2mm 17mm 9.5mm}, clip]{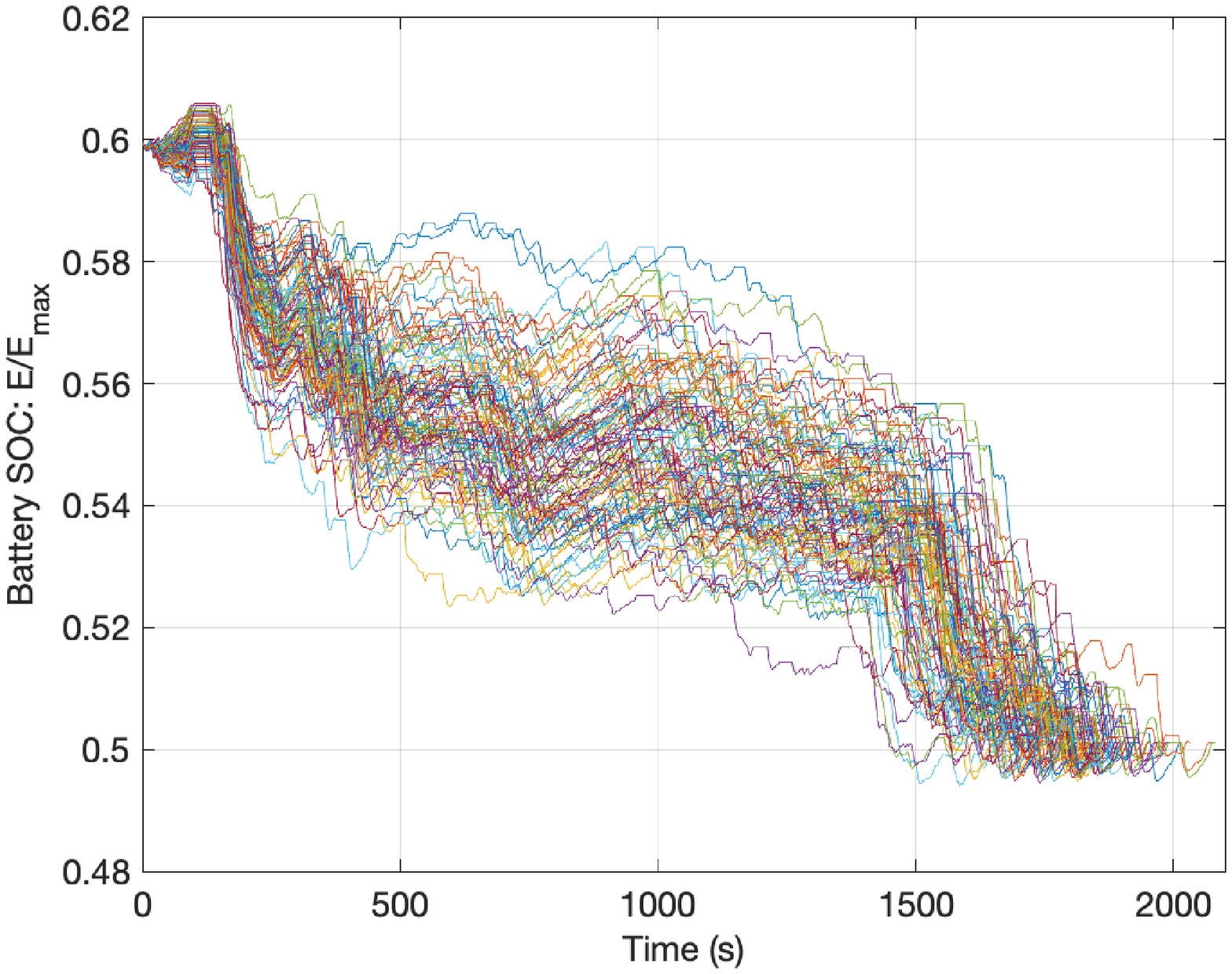}}
\vspace{-2mm}
\caption{Predicted battery state of charge profiles corresponding to the optimal solution of~(\ref{eq:PHEV_opt}) with $100$ power demand scenarios generated from random perturbations of the FTP-75 drive cycle.}
\label{fig:ADMM100}
\vspace{-2mm}
\end{figure}

Thresholds for the primal and dual residuals were set at $\bar{r} = 10^{-6} \Delta E$ and $\bar{\sigma} = 10^{-2}$. 
%
%
Step-length parameters $\rho_1$, $\rho_2$, $\rho_3$, $\rho_4$ were adapted at each ADMM iteration according to
\[
\rho_i \gets \tau \rho_i \ \text{if} \ r/\sigma > 1.2 \bar{r}/\bar{\sigma} 
\ \ \text{and} \ \ 
\rho_i \gets \rho_i/\tau \ \text{if} \  r/\sigma < 0.8 \bar{r}/\bar{\sigma}
\]
with $\tau = 1.1$ and initial values $\rho_1 = 10^{-4}$, $\rho_2 = 2\times 10^{-6}$, $\rho_3 = \rho_4 = 5\times 10^{-6}$. 
This adaptive approach, which is based on the discussion in Sec.~3.4.1 of \cite{boyd11}, was found to reduce the effect of  the value of $q$ on the rate of convergence of the ADMM iteration.
A typical solution for the battery state of charge for $q=100$ demand scenarios is shown in Fig.~\ref{fig:ADMM100}.

\subsection{Comparison of CPU/GPU performance}
\label{sec:comp_admm}

\noindent
The CPU-based ADMM implementation was tested using the following procedure:
\begin{enumerate}
  \item Perform steps~1 and 2 of Algorithm~\ref{alg:PHEV_control}.
  \item Record start time $T_1$.
  \item Execute step 3 of  Algorithm~\ref{alg:PHEV_control} on the CPU.
  \item Record end time $T_2$ and elapsed time $T$ = $T_2 - T_1$.
\end{enumerate}
For the parallel implementation, step 3 was replaced by:
\begin{enumerate}
  \setcounter{enumi}{2}
  \item
  \begin{enumerate}
  	\item Input data to the GPU.
  	\item Execute step 3 of Algorithm~\ref{alg:PHEV_control} on the GPU.
  	\item Copy results back to the CPU.
  \end{enumerate}
\end{enumerate}

\begin{figure}[t]
\vspace{-2mm}
\centerline{\includegraphics[scale = 0.47, trim={8mm 0mm 2mm 9mm}, clip]{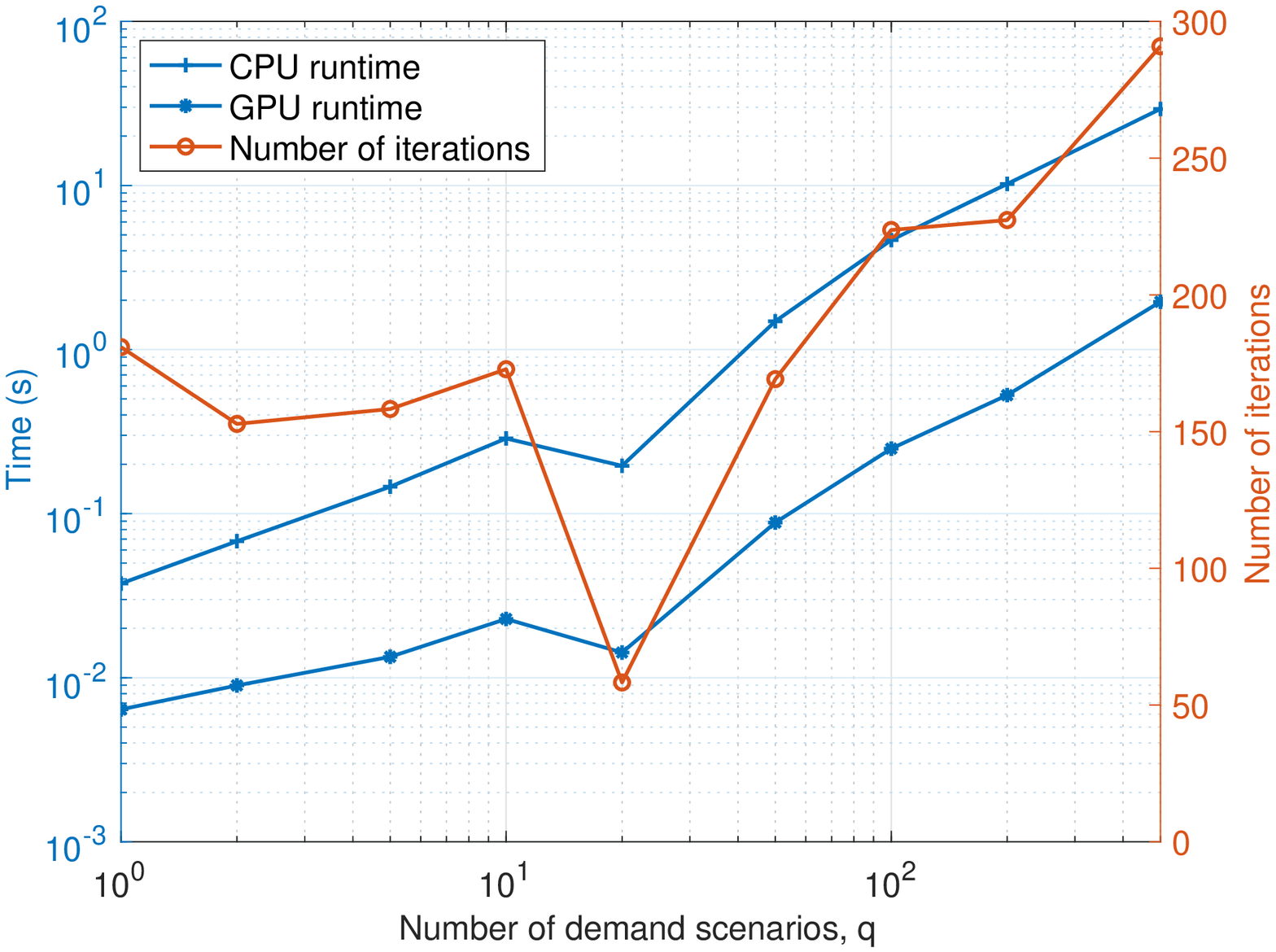}}
\vspace{-2mm}
\caption{Mean runtimes of CPU and GPU (parallel) implementations (left vertical axis) and number of ADMM iterations (right vertical axis) against number of demand scenarios (horizontal axis).}\label{fig:admm_tests}
\vspace{-2mm}
\end{figure}

For numbers of demand scenarios $q$ between 5 and 500, Figure \ref{fig:admm_tests} shows that GPU runtimes were between 10 and 20 times faster than tests executed on a 4-core CPU with (/Ox) speed optimizations.
%
The number of ADMM iterations required tends to increase with $q$, although for the chosen parameters a significant reduction is observed around $q=20$. 
The speed-up achieved by the GPU implementation increases with $q$, even though for $q>10$ the number of threads required to execute the ADMM iteration exceeds the limit that can be executed simultaneously on the GPU.
To reduce the effect on overall performance of computing primal and dual residuals, the termination conditions were checked only once every 10 iterations in both implementations.  
%

\section{Conclusions}
This paper describes an ADMM algorithm for robust quadratic resource allocation problems. We discuss a parallel implementation of the algorithm and its application to optimal energy management for PHEVs with uncertain power demand.   
%
Our results suggest that it is feasible to solve optimal PHEV energy management problems online using inexpensive computing hardware.
We also note that the state of the art in low-cost parallel processing hardware is evolving rapidly and that next generation GPUs can provide significant speed improvements over the results reported here. 

\section*{Appendix A: ADMM iteration}

To express~(\ref{eq:robust_opt}) in a form suitable for ADMM, we introduce slack variables $s_k^{(j)}$ into the demand constraints and introduce dummy variables $\smash{z_k^{(i,j)}}$, $\smash{h^{(i,j)}}$ into the capacity constraints: \begin{equation}\label{eq:admm_opt}
\begin{aligned}
& \min_{\substack{x_1^{(i)}\in\mathbb{R}, \ x_k^{(i,j)}\in \lubx, \ z_k^{(i,j)}\in\mathbb{R}\\ s_k^{(j)}\geq 0, \ h^{(i,j)} \leq c^{(i)}}} \ \ 
\tfrac{1}{q} \sum_{i,j,k} f _k^{(i,j)}(x_k^{(i,j)}) \\
& \text{subject to}\ 
\begin{alignedat}[t]{2}
& \sum_{i} x_k^{(i,j)} - y_k^{(j)}  = s_k^{(j)}
& j \in\N_q, \, k\in\N_n \\
& g_k^{(i,j)}(x_k^{(i,j)}) = z_{k}^{(i,j)} 
& i\in\N_n, \, j \in\N_q, \, k\in\N_n \\
& \sum_{k} z_k^{(i,j)} = h^{(i,j)} 
& i\in\N_n, \, j \in\N_q \\
& x_1^{(i,j)} = x_1^{(i)} 
& i\in\N_n, \, j \in \N_q
\end{alignedat}
\end{aligned}
\end{equation}
This formulation preserves separability of the nonlinear cost and constraints by imposing linear constraints on ${z_k^{(i,j)}}$.
Let $\ind_{\mathbb{I}}(\cdot)$ denote the indicator function of the interval $\mathbb{I}$ 
(for scalar $x$, $\ind_{\mathbb{I}}(x) = 0$ if $x\in \mathbb{I}$ and $\ind_{\mathbb{I}}(x) = \infty$ if $x\notin \mathbb{I}$; for vector $x$, $\ind_{\mathbb{I}}(x)$ is the vector with $k$th element equal to $\ind_{\mathbb{I}}(x_k)$). We define the augmented Lagrangian function 
\begin{align*}
& L(x_1,x,z,s,h,\mu,\lambda,\nu,p) := 
\tfrac{1}{q} \sum_{i,j,k} f_k^{(i,j)}(x_k^{(i,j)}) \\
& +
\sum_{i,j,k} \Bigl[ 
\ind_{[\underline{x}_k^{(i)}\hspace{-0.5ex},\hspace{0.3ex}\bar{x}_k^{(i)}]}
(x_k^{(i,j)}) \!+\! \tfrac{\rho_1}{2} \bigl( z_k^{(i,j)} \!-\! g_k^{(i,j)} (x_k^{(i,j)}) \!+\! \lambda_k^{(i,j)} \bigr)^2 \Bigr]
\\ 
& + 
\sum_{i,j} \Bigl[ 
\ind_{(-\infty,\hspace{0.3ex}c^{(i)}]}(h^{(i,j)}) + \tfrac{\rho_2}{2} \bigl( h^{(i,j)} - \ones^{\top} z^{(i,j)} + p^{(i,j)}\bigr)^2\Bigr]
\\
& +
\sum_{j} \Bigl[ \ones^\top \ind_{[0,\infty)}(s^{(j)}) + \tfrac{\rho_3}{2} \bigl\| s^{(j)} \!-\! \sum_{i} x^{(i,j)} \!+\! y^{(j)} \!+\! \mu^{(j)} \bigr\|^2 \Bigr]
\\
& +
{\sum_{i,j}} \tfrac{\rho_4}{2} \bigl( x_1^{(i)} - x_1^{(i,j)} + \nu^{(i,j)} \bigr)^2
\end{align*}
Here 
$\mu^{(j)}$,
$\lambda^{(i,j)}$, 
$\nu^{(i,j)}$, 
$p^{(i,j)}$ and
$\nu^{(i,j)}$ are multipliers for the equality constraints 
in (\ref{eq:admm_opt});
$\smash{y^{(j)}}$, $\smash{x^{(i,j)}}$, $\smash{z^{(i,j)}}$, $\smash{s^{(j)}}$ and $\smash{g^{(i,j)}(x^{(i,j)})}$ are the vectors with $k$th elements equal to $\smash{y_k^{(i)}}$, $\smash{x_k^{(i,j)}}$, $\smash{z_k^{(i,j)}}$, $\smash{s_k^{(j)}}$ and $\smash{g_k^{(i,j)}(x_k^{(i,j)})}$, and $\ones = [\, 1 \ \cdots \ 1\,]^\top$. 
Also $\rho_1, \ldots \rho_4$ are non-negative scalars that control the step length of multiplier updates in the ADMM iteration~\cite[Sec.~3.1]{boyd11}. These design parameters control the rate of convergence and can be tuned to the problem at hand  (see e.g.~\cite{east18}).

At each iteration, $L$ is minimized with respect to primal variables 
$x^{(i,j)}$, $x_1^{(i)}$, $z^{(i,j)}$, $s^{(j)}$, $h^{(i,j)}$. The primal and dual variable updates are given, for $k\in\N_n$, $i\in\N_m$, $j\in\N_q$, by
\begin{subequations}\label{eq:admm_iter}
\begin{align}
& x_k^{(i,j)} \gets \Pi_{\lubx} \Bigl( \begin{aligned}[t]  
&  \argmin_{x^{(i,j)}_k} \Bigl\{  \tfrac{1}{q} f_k^{(i,j)}(x_k^{(i,j)}) \\
& +
\tfrac{\rho_1}{2} \bigl( z_k^{(i,j)} - g_k^{(i,j)}(x_k^{(i,j)}) + \lambda_k^{(i,j)}\bigr)^2 \\
& + \tfrac{\rho_3}{2} \bigl( s_k^{(j)} - \sum_{i} x_k^{(i,j)}  + y_k^{(j)} + \mu_k^{(j)}\bigr)^2 \\
& + \delta_{k,1}\tfrac{\rho_4}{2} \bigl( x_1^{(i)} - x_1^{(i,j)} + \nu^{(i,j)}\bigr)^2\Bigr\}\Bigr)
\end{aligned}
\\
& \begin{aligned}
& z^{(i,j)} \gets g^{(i,j)}(x^{(i,j)}) - \lambda^{(i,j)} \\
&\hspace{2.5em} +\!\frac{\rho_2}{\rho_1 \!+\! n\rho_2} \ones \Bigl(h^{(i,j)} \!+\! p^{(i,j)}  \!-\! \ones^{\!\top}\!\bigl(g^{(i,j)}(x^{(i,j)}) \!-\! \lambda^{(i,j)}\bigr)\Bigr)
\end{aligned}
\\
& x_1^{(i)} \gets \sum_j \bigl( x_1^{(i,j)} - \nu^{(i,j)} \bigr)
\\
& h^{(i,j)} \gets \Pi_{(-\infty, c^{(i)}]}\bigl(\ones^\top z^{(i,j)} - p^{(i,j)} \bigr)
\\
& s^{(j)} \gets \Pi_{[0,\infty)} \Bigl(\, \sum_{i=1}^m x^{(i,j)} - y^{(j)} - \mu^{(j)}\Bigr)
\\
& \mu^{(j)} \gets \mu^{(j)} + s^{(j)} - \sum_{i} x^{(i,j)} + y^{(j)}
\\
& \lambda^{(i,j)} \gets \lambda^{(i,j)} + z^{(i,j)} - g^{(i,j)}(x^{(i,j)})
\\
& \nu^{(i,j)} \gets \nu^{(i,j)} + x_1^{(i)} - x_1^{(i,j)}
\\
& p^{(i,j)} \gets p^{(i,j)} + h^{(i,j)} - \ones^\top z^{(i,j)} 
\end{align}
\end{subequations}
where $\Pi_{\mathbb{I}}(x)$ for scalar $x$ denotes projection onto the interval $\mathbb{I}$, while for vector $x$, $\Pi_{\mathbb{I}}(x)$ is the vector obtained by  projecting each element of $x$ onto $\mathbb{I}$, and $\delta_{k,1} = 1$ if $k=1$, $\delta_{k,1} = 0$ if $k\neq 1$.
Note that the update (\ref{eq:admm_iter}a) can be written equivalently as
\begin{align*} 
\!\!\! x_k^{(i,j)} &\gets \Pi_{\lubx} \Bigl(  \argmin_{x^{(i,j)}_k} \Bigl\{ 
\tfrac{1}{q} f_k^{(i,j)}(x_k^{(i,j)})
\\
& + \tfrac{\rho_1}{2} \bigl( \theta_k^{(i,j)} \!-\! g_k^{(i,j)}(x_k^{(i,j)}) \bigr)^2 
+ \tfrac{\rho_3}{2} \big( \phi_k^{(i,j)} \!-\! x_k^{(i,j)}\bigr)^2 
\\
&\hspace{5em} + \delta_{k,1}\tfrac{\rho_4}{2} \bigl( x_1^{(i)} - x_1^{(i,j)} + \nu^{(i,j)}\bigr)^2\Bigr\}\Bigr) .
\end{align*}
with $\theta_k^{(i,j)} = z_k^{(i,j)} + \lambda_k^{(i)}$ and
$\phi_k^{(i,j)} = s_k^{(j)} - \sum_{l\neq i} x_k^{(l,j)} + y_k^{(j)} + \mu_k^{(j)}$.
The iteration terminates when the primal and dual residuals $r$ and $\sigma$ defined by
\begin{align*}
r &= \max \begin{aligned}[t] \Bigl\{ 
& \max_j \, \bigl\|s^{(j)}  - \textstyle\sum_{i} x^{(i,j)} + y^{(j)} \bigr\|_\infty, \\
& \max_{i,j} \, \bigl\|z^{(i,j)} - g^{(i,j)}(x^{(i,j)}) \bigr\|_\infty , \\
& \max_{i,j} \bigl\lvert h^{(i,j)} - \ones^\top z^{(i,j)} \bigr\rvert ,
\max_{i,j} \, \bigl\lvert x_1^{(i,j)} - x_1^{(i)} \bigl\rvert \Bigr\} \end{aligned} 
\\
\sigma &=
\max \begin{aligned}[t] \Bigl\{
& \rho_1 \max_i \bigl\| z^{(i)} - \tilde{z}^{(i)} \bigr\|_\infty  ,
\ \rho_2 \max_{i,j}\bigl\lvert h^{(i,j)} - \tilde{h}^{(i,j)} \bigl\rvert ,
\\
& \rho_3 \max_{i,j} \bigl\| s^{(j)} -  \tilde{s}^{(j)} - \textstyle\sum_{i} (x^{(i,j)} - \tilde{x}^{(i,j)} ) \bigr\|_\infty \Bigr\} \end{aligned}
\end{align*}
(where $[(\cdot)]_0$ denotes the value of an iterate $(\cdot)$ at the previous iteration), fall below predefined thresholds, $\bar{r}$ and $\bar{\sigma}$.

\bibliographystyle{IEEEtran}
\bibliography{refs_new}

\end{document}